\documentstyle[12pt]{article}
\title{ THE STRONGLY SYMMETRIC ELEMENTS AND  SOLUTIONS OF
 YANG-BAXTER EQUATION
 \thanks {This work is supported by National Science Foundation}}
\author{Shouchuan Zhang \\ Department  of Mathematics, Nanjing  University,
210008 \\ P.R.China}
\date{}
\begin{document}
\newtheorem{Theorem}{\quad Theorem}[section]
\newtheorem{Proposition}[Theorem]{\quad Proposition}
\newtheorem{Definition}[Theorem]{\quad Definition}
\newtheorem{Corollary}[Theorem]{\quad Corollary}
\newtheorem{Lemma}[Theorem]{\quad Lemma}
\newtheorem{Example}[Theorem]{\quad Example}
\maketitle
\begin {abstract}

 It is shown that all strongly symmetric  elements are solutions of
  constant classical Yang-Baxter equation  in Lie algebra,     or               of
  quantum  Yang-Baxter  equation  in algebra.
Otherwise, all  solutions of constant classical Yang-Baxter
 equation (CYBE)  in Lie algebra $L$ with dim $L \le 3$ over field $k$ of
 characteristic $2$  are obtained
\end {abstract}
\noindent
\addtocounter{section}{-1}

\section{ Introduction}

The Yang-Baxter equation first came up in a paper by Yang as
factorition condition of the scattering S-matrix in
the many-body  problem in one dimension and in work of Baxter on
exactly solvable models in statistical mechanics.
It has been playing an important role in mathematics and physics
( see \cite {BD} , \cite {YG} ).
Attempts to find solutions of The Yang-Baxter equation in a systematic
way have led to the theory of quantum groups.
The Yang-Baxter
equation is of many forms.
The classical Yang-Baxter equation  and quantum Yang-Baxter equation are two
kinds of them.

 In many applications one need to know the solutions of the two equations.
 The author in \cite {Z} systematically studied the solutions of low
dimensional Lie algebras, found strongly symmetric elements and shew
that all
of them are solutions of CYBE. Naturally,
we would like to ask
the    following questions:

(1) Is every strongly element in $L \otimes L$  a solution of CYBE for
any dimensional Lie algebra $L$?

(2) Do the conclusions in (BI2494R)  hold for field of characteristic 2 ?

(3) Is every strongly symmetric element  in $A \otimes A$ a solution of
QYBE for any algebra $A$?

The resolution of these questions is very necessary for not only physics
but also mathematics.

To answer these questions,
the author  write the paper.
 In this paper,
the author  shows  that all strongly symmetric  elements are solutions of
  constant classical Yang-Baxter equation  in Lie algebra,     or               of
  quantum  Yang-Baxter  equation  in algebra.
   Otherwise, to bring the paper \cite {Z} completion, the author find
     all  solutions of CYBE  in Lie algebra $L$ with dim $L \le 3$ over field $k$ of
 characteristic $2$.

Let $k$ be a field and $L$ be a Lie algebra over $k$.
If $r = \sum a_i \otimes b_i  \in L \otimes L $ and $x\in L$, we define
 \begin{eqnarray*}                              
\left[ r^{12}, r^{13} \right] &:=& \sum_{i, j} [a_i, a_j] \otimes b_i
\otimes b_j \\
\left[ r^{12}, r^{23} \right] &:=& \sum_{i, j} a_i\otimes [b_i,a_j]
 \otimes  b_j \\
\left[ r^{13}, r^{23} \right] &:=& \sum_{i, j} a_i\otimes a_j \otimes [b_i, b_j] \\
\end{eqnarray*}

and call $$     [r^{12}, r^{13}]+  [r^{12}, r^{23}]+ [r^{13}, r^{23}]=0$$
the classical Yang-Baxter equation (CYBE).

Let $A$ be an algebra and  $R = \sum a_i \otimes b_i  \in A \otimes  A$.
We define
 \begin{eqnarray*}                              
R^{12} &=& \sum _i a_i \otimes b_i \otimes 1  \\
R^{13} &=& \sum _i a_i \otimes 1\otimes b_i   \\
R^{23} &=& \sum _i 1 \otimes  a_i \otimes b_i
\end{eqnarray*}

and call $$
R^{12}R^{13}R^{23} = R^{23}R^{13}R^{12} $$
the quantum  Yang-Baxter equation (QYBE).

\section {Strongly symmetric elements  }
In this section,                             
we show  that all strongly symmetric  elements are solutions of
  constant classical Yang-Baxter equation  in Lie algebra,     or               of
  quantum  Yang-Baxter  equation  in algebra.

\begin {Definition} \label {0.1}
 Let  $\{ e_i \mid i \in \Omega \}$ be a basis of vector space $V$
and     $r = \sum_{i,j \in \Omega }
 k_{ij}(e_{i} \otimes e_{j}) \in V \otimes V,$
where   $k_{ij} \in k,$

If  $$k_{ij}=k_{ji}  \hbox { \ \ \ } k_{ij} k_{lm} = k_{il}k_{jm}$$
for $i, j, l, m \in \Omega$
then  $r$ is called
strongly symmetric to the basis
 $\{ e_i \mid  i \in \Omega  \}$
                  \end {Definition}

Obviously,  symmetry does not depend on the particular choice of basis
of $V$ . We need to know if
strong  symmetry  depends on the particular choice of basis of $V$.

\begin {Lemma}                                            
\label {0.2}
 Let  $\{ e_i \mid i \in \Omega \}$ be a basis of vector space $V$
and     $r = \sum_{i,j \in \Omega } k_{ij}(e_{i} \otimes e_{j}) \in V \otimes V,$
where   $k_{ij} \in k.$  Then

(I) The strong symmetry  does not depend on the particular choice of
basis of $V;$

(II)  $r$  is strongly symmetric iff
   $$   k_{ij} k_{lm} = k_{il}k_{jm}$$
for  any  $i, j, l, m \in \Omega $

(III)  $r$ is strongly symmetric iff
   $$   k_{ij} k_{lm} = k_{st}k_{uv}$$
where $\{ s, t, u, v \}$ is a  combination of  $\{i, j, l, m \} $
for  any  $i, j, l, m \in \Omega $

(IV)  $r$ is strongly symmetric iff
   
$r=0;$  or

$r = \sum_{i,j \in \Omega }  k_{i_0 i}k_{i_0j} k_{i_0i_0}^{-1}
        (e_{i} \otimes e_{j}) \in V \otimes V,$
where  $k_{i_0 i_0} \not=0;$

\end {Lemma}

{ \bf Proof. }
(I)  Let  $\{  e_i' \mid i \in \Omega ' \}$ be another  basis of $V$
It is sufficient to show that  $r$ is strongly symmetric to
the basis  $\{  e_i' \mid i \in \Omega ' \}$.
Obviously, there exists $q_{ij} \in k$
such that  $e_{i} = \sum _{s} e_s'q_{si}  $    for $i \in \Omega.$
By computation,  we have that
$$r = \sum_{ s, t  \in \Omega }( \sum _{i,j \in \Omega '} k_{ij} q_{si}q_{tj} )(e_{s}' \otimes e_{t}').$$
If set   $k_{st}'= \sum _{i,j} k_{ij} q_{si}q_{tj}$,
then
for $l, m, u, v \in \Omega '  $   we have that
 \begin {eqnarray*}
 k_{lm}'k_{uv}'
  &=&  \sum_{i,j,s,t} k_{ij} k_{st} q_{li}q_{mj} q_{us}q_{vt}   \\
  &=&  \sum_{i,j,s,t} k_{is} k_{jt} q_{li}q_{mj} q_{us}q_{vt} { \ \ }
 (  \hbox { \ \  by  }   k_{ij}k_{st} = k_{is}k_{jt}  ) \\
&=&  k_{lu}'k_{mv}'
     \end {eqnarray*}
 and
 \begin {eqnarray*}
 k_{lm}'&=&k_{ml}'
 \end {eqnarray*}
which implies that
 $r$ is strongly symmetric to
the basis $\{  e_i' \mid i \in \Omega ' \}$
                                                  
(II) The necessity is obvious. Now we sufficiency. For any $i, j
\in \Omega  $, if $k_{ii} \not=0$, then $k_{ij} = k_{ji}$
since $ k_{ij}k_{ii} = k_{ii}k_{ji}$. If $k_{ii} = 0$, then
$k_{ij}= k_{ji} = 0$ since $k_{ij} ^2 = k_{ii}k_{jj} =0 = k_{ji}^2.$
                                                  
(III) The sufficiency follows from part (II). Now we show the necessity.
 For any $i, j , l, m
\in \Omega  $, we see that
 $k_{ij}k_{lm} = k_{ij}k{ml} = k_{il}k{jm} =k_{il}k{mj} =k_{im}k{lj}
 =k_{im}k{lj} $

 Similarly, we can get that $k_{ij}k_{lm}=k_{st}k_{uv}$  when $s=j$ or $s=l$
 or $s=m$.

(IV)  If
$r$ is a case in part (IV),
 then $r$  is strongly symmetric by straightforward verification. Conversely,
 if $r$ is strongly symmetric and $r\not=0$, then there exists
 $i_0 \in \Omega$  such that  $k_{i_0i_0} \not=0.$
 Thus $k_{ij}  =k_{i_0i}k_{i_0j} k_{i_0i_0}^{-1}$
$\Box$

If $A$ is an algebra then we can get a Lie algebra $L(A)$ by defining
$[x,y] =xy -yx $ for any $x, y \in A$.

\begin {Theorem} \label {0.3}      Let $A$  be an algebra and $r \in A\otimes A$ be strongly symmetric.
Then  $r$ is a solution of CYBE in $L(A)$ and of QYBE in $A$.
\end {Theorem}
{\bf Proof.}   Let $\{ e_i \mid i \in \Omega \}$  be a basis of $A$ and
the multiplication of $A$ be defined as follows:
$$e_i e_j = \sum _m a_{ij}^m e_m $$
for any $i, j \in \Omega $

  We first show that $r$ is  a solution of CYBE in Lie algebra $L(A)$.
   By computation, for all $i, j, l \in \Omega ,$   we have that the coefficient of
  $e_i \otimes e_j \otimes e_l$ in
$[r^{12}, r^{13}] $
is               
\begin {eqnarray*}
&{}&\sum _{s t} a_{s,t} ^ik_{sj}k_{tl}-a_{ts} ^ik_{sj}k_{tl}   \\
&=& \sum _{s, t} a_{st}^i k_{sj}k_{tl}-
\sum _{s,t} a_{st}^i k_{tj}k_{sl}    \\
&=& 0   \hbox { \ \ \ \ \ \  by Lemma \ref {0.2} (III)}
\end {eqnarray*}
Similarly,  we have that the coefficients of
  $e_i \otimes e_j \otimes e_l$ in
$[r^{12}, r^{23}] $ and
 $[r^{13}, r^{23}] $ are zero.
Thus $r$ is  a solution  of CYBE in $L(A)$.

  Next we show that $r$ is a solution of QYBE in algebra $A$.
We denote $r$ by $R$ since  the solution of QYBE is usually donoted by $R$.

 By computation, for all $i, j, l \in \Omega ,$   we have that the coefficient of
  $e_i \otimes e_j \otimes e_l$ in
$R^{12}R^{13} R^{23} $
is               
\begin {eqnarray*}
&{}&\sum _{s, t, u,v, w,m}
a_{st} ^i
a_{um} ^j
a_{vw} ^l k_{su}k_{tv}k_{mw}
\end {eqnarray*}

and the coefficient of
  $e_i \otimes e_j \otimes e_l$ in
$R^{23}R^{13} R^{12} $
is               
\begin {eqnarray*}
&{}&\sum _{s, t, u,v, w,m}
a_{tm} ^i
a_{sw} ^j
a_{uv} ^l k_{su}k_{tv}k_{mw} \\
&=&
\sum _{s, t, u,v, w,m}
a_{st} ^i
a_{um} ^j
a_{vw} ^l k_{su}k_{tv}k_{mw}  \hbox { \ \ \  by Lemma \ref {0.2} (III)}
\end {eqnarray*}
 Thus $r$ is a solution of QYBE in algebra $A.$  $\Box$

\begin {Corollary}  \label {0.4}   Let $L$ be a Lie algebra and $r\in
L\otimes L$
be strongly symmetric. Then
$r$ is a solution of CYBE.
\end {Corollary}

{\bf Proof:}
Let $A$ denote the unversal enveloping algebra of $L$ . It is clear
that  $r$ is also strongly symmetric in $A \otimes A$.
Thus

 $$     [r^{12}, r^{13}]+  [r^{12}, r^{23}]+ [r^{13}, r^{23}]=0$$
by Theorem \ref {0.3} . Consequently, $r$ is a solution of
CYBE in $L$. $\Box$

\section {The solutions of CYBE }
In this section, we find the general solution of CYBE for Lie algebra $L$
over field $k$ of characteristic 2
with $dim \ \ L \le 3.$ Throughtout  this section, the characteristic of $k$
is 2.

\begin {Proposition}\label {1.3}
Let $L$ be a Lie algebra with $dim { \ } L =2$.
Then $r$ is a  solution of CYBE iff
$r$ is  symmetric.
\end {Proposition}

{\bf Proof.}   It is not hard since the tensor $r$ has only 4 coefficients.
 $\Box$

  If $r = x (e _1 \otimes e_1) + y (e _2 \otimes e_2) + 
z (e _3 \otimes e_3) + p (e _1 \otimes e_2) +p(e _2 \otimes e_1)+
s (e _1 \otimes e_3)+s(e _3 \otimes e_1)+ u (e _2 \otimes e_3)
+u (e _3 \otimes e_2)$  with  $\alpha y z + \beta xz + xy +
  \beta s^2 + \alpha u^2 + p^2 =0,$
then  $r$ is called
 $\alpha , \beta $-symmetric to the
basis  $\{ e_1, e_2, e_3 \}.$

\begin {Proposition}
\label {1.4} Let
$L$ be a Lie algebra with a basis $\{ e_1, e_2, e_3 \}$
such that
$[e_1,e_2]= e_3, [e_2,e_3]=\alpha e_1, [e_3,e_1]=\beta e_2$, where
$\alpha, \beta \in k$. Let $ p, q, s, t, u, v, x, y, z \in k$.

Then  (I) $r$ is a solution of CYBE in $L$ iff $r$ is
 $\alpha, \beta $- symmetric
to the basis $\{ e_1, e_2, e_3 \}.$

(II) If $r$ is strongly symmetric, then $r$ is a solution of CYBE in $L$
\end {Proposition}

{ \bf Proof .}
Let   $r = \sum_{i,j= 1 }^{3} k_{ij}(e_{i} \otimes e_{j}) \in L \otimes L$
and  $k_{ij} \in k,$ with  $i, j =1, 2, 3.$
It is clear that
\begin {eqnarray*}
\left[ r^{12}, r^{13} \right] &=& \sum _{i,j=1}^3 \sum _{s,t =1}^3 k_{ij}k_{st}
 [e_i,e_s] \otimes e_j \otimes e_t  \\
\left[r^{12}, r^{23} \right] &=& \sum _{i,j=1}^3 \sum _{s,t =1}^3 k_{ij}k_{st}
e_i \otimes[ e_j, e_s ]\otimes e_t \\
\left[r^{13}, r^{23} \right] &=& \sum _{i,j=1}^3 \sum _{s,t =1}^3 k_{ij}k_{st}
e_i \otimes e_s \otimes[e_j, e_t]
\end {eqnarray*}
 By computation, for all $i, j, n =1, 2, 3,$   we have that the coefficient of
  $e_j \otimes e_i \otimes e_i$ in
$[r^{12}, r^{13}] $
is zero and
 $e_i \otimes e_i \otimes e_j$
in  $[r^{13}, r^{23}]$ is zero.

We now see the coefficient
 of  $e_i \otimes e_j \otimes e_n$

   in  $[r^{12}, r^{13}] + [r^{12}, r^{23}] +[r^{13}, r^{23}]$.

(1) $e_1 \otimes e_1 \otimes e_1 (
\alpha k_{12}  k_{31}-
\alpha k_{13} k_{21})$

(2) $e_2 \otimes e_2 \otimes e_2 (
 - \beta k_{21}  k_{32}+\beta k_{23} k_{12} )$

(3) $e_3 \otimes e_3 \otimes e_3 (
k_{31}  k_{23}- k_{32} k_{13})$

(4) $e_1 \otimes e_2 \otimes e_3 (
\alpha k_{22}k_{33} - \alpha k_{32} k_{23} -  \beta k_{11}  k_{33}
+\beta k_{13} k_{13}
+ k_{11}  k_{22} -k_{12} k_{21} )$

(5) $e_2 \otimes e_3 \otimes e_1 (
\beta k_{33}k_{11} - \beta k_{13} k_{31} - k_{22}  k_{11}+ k_{21} k_{21}
+ \alpha k_{22}  k_{33} -\alpha k_{23} k_{32} )$

(6) $e_3 \otimes e_1 \otimes e_2 (
k_{11}k_{22} - k_{21} k_{12} - \alpha k_{33}  k_{22}+ \alpha k_{32} k_{32}
+ \beta k_{33}  k_{11} -\beta k_{31} k_{13} )$

(7) $e_1 \otimes e_3 \otimes e_2 (
-\alpha k_{33}k_{22} + \alpha k_{23} k_{32} + k_{11}  k_{22}- k_{12} k_{12}
-\beta k_{11}  k_{33} + \beta k_{13} k_{31} )$

(8) $e_3 \otimes e_2 \otimes e_1 (
-k_{22}k_{11} + k_{12} k_{21} + \beta k_{33}  k_{11}- \beta k_{31} k_{31}
- \alpha k_{33}  k_{22} + \alpha k_{32} k_{23} )$

(9) $e_2 \otimes e_1 \otimes e_3 (
-\beta k_{11}k_{33} + \beta k_{31} k_{13}  + \alpha k_{22}  k_{33}
- \alpha k_{23} k_{23}
- k_{22}  k_{11} +k_{21} k_{12} )$

(10) $e_1 \otimes e_1 \otimes e_2 (
- \alpha k_{31}k_{22} + \alpha k_{21} k_{32} + \alpha k_{12}  k_{32}
- \alpha k_{13} k_{22})$

(11) $e_2 \otimes e_1 \otimes e_1 (
 -\alpha k_{23}  k_{21}+
\alpha k_{22} k_{31} +\alpha k_{22}  k_{13} - \alpha k_{23} k_{12} )$

(12) $e_1 \otimes e_1 \otimes e_3 (
\alpha k_{21}k_{33} - \alpha k_{31} k_{23} + \alpha k_{12}  k_{33}
 - \alpha k_{13} k_{23} )$

(13) $e_3 \otimes e_1 \otimes e_1 (
 \alpha k_{32}  k_{31}- \alpha k_{33} k_{21}
+ \alpha k_{32}  k_{13} -\alpha k_{33} k_{12} )$

(14) $e_2 \otimes e_2 \otimes e_1 (
- \beta k_{12}k_{31} + \beta k_{32} k_{11} + \beta k_{23}  k_{11}
-\beta  k_{21} k_{31})$

(15) $e_1 \otimes e_2 \otimes e_2 (
 -\beta  k_{11}  k_{32}+
\beta k_{13} k_{12} -\beta k_{11}  k_{23} + \beta k_{13} k_{21} )$

(16) $e_2 \otimes e_2 \otimes e_3 (
\beta k_{32}k_{13} -\beta k_{12} k_{33} -\beta k_{21}  k_{33}
+\beta k_{23} k_{13})$

(17) $e_3 \otimes e_2 \otimes e_2 (
-\beta k_{31}  k_{32}+ \beta k_{33} k_{12}
+ \beta k_{33}  k_{21} - \beta k_{31} k_{23} )$

(18) $e_3 \otimes e_3 \otimes e_1 (
k_{13}k_{21} - k_{23} k_{11} + k_{31}  k_{21}- k_{32} k_{11} )$

(19) $e_3 \otimes e_3 \otimes e_2 (
-k_{23}k_{12} + k_{13} k_{22} - k_{32}  k_{12}+ k_{31} k_{22})$

(20) $e_2 \otimes e_3 \otimes e_3 (
 k_{21}  k_{23} - k_{22} k_{13}
+ k_{21}  k_{32} -k_{22} k_{31} )$

(21) $e_1 \otimes e_3 \otimes e_3 (
- k_{12}  k_{13}+ k_{11} k_{23}
- k_{12}  k_{31} +k_{11} k_{32} )$

(22) $e_1 \otimes e_3 \otimes e_1 (
\alpha k_{23}k_{31} - \alpha k_{33} k_{21} - k_{12}  k_{11}+ k_{11} k_{21}
+\alpha k_{12}  k_{33} -\alpha k_{13} k_{32} )$

(23) $e_1 \otimes e_2 \otimes e_1 (
-\alpha k_{32}k_{21} +\alpha  k_{22} k_{31} - \beta k_{11}  k_{31}
+ \beta k_{13} k_{11} -\alpha k_{13}  k_{22} +\alpha k_{12} k_{23} )$

(24) $e_2 \otimes e_1 \otimes e_2 (
\beta k_{31}k_{12} -\beta k_{11} k_{32} + \alpha k_{22}  k_{32}
-\alpha  k_{23} k_{22} - \beta k_{21}  k_{13} +\beta k_{23} k_{11} )$

(25) $e_2 \otimes e_3 \otimes e_2 (
-\beta k_{13}k_{32} + \beta k_{33} k_{12} + k_{21}  k_{22}- k_{22} k_{12}
-\beta  k_{21}  k_{33} +\beta k_{23} k_{31} )$

(26) $e_3 \otimes e_2 \otimes e_3 (
k_{12}k_{23} - k_{22} k_{13} -\beta k_{31}  k_{33}+ \beta k_{33} k_{13}
+ k_{31}  k_{22} -k_{32} k_{21} )$

(27) $e_3 \otimes e_1 \otimes e_3 (
-k_{21}k_{13} + k_{11} k_{23} + \alpha k_{32}  k_{33}- \alpha k_{33} k_{23}
+ k_{31}  k_{12} -k_{32} k_{11} ).$
    
Let $k_{11}=x, k_{22}=y,  k_{33}=z, k_{12}=p,  k_{21} =q, k_{13}=s,
 k_{31} =t, k_{23}=u,  k_{32}=v.$

It follows from (1)-(27)  that

(28) $\alpha pt= \alpha qs $

(29) $\beta qv= \beta pu$

(30) $tu=vs$

(31) $  \alpha yz - \beta xz + xy - \alpha uv + \beta s^2 - pq =0.$

(32) $ \beta zx -yx+ \alpha yz - \beta st + q^2-\alpha uv =0$

(33) $ xy  -\alpha zy + \beta zx -  pq+ \alpha v^2 - \beta st =0$

(34) $- \alpha zy+ xy- \beta xz+  \alpha uv - p^2+ \beta st =0$

(35) $-xy +\beta xz- \alpha yz + pq - \beta t^2+ \alpha uv =0$

(36) $ -\beta xz+ \alpha yz - yx+  \beta st -  \alpha u^2 + pq =0$

(37)  $\alpha (-ty + qv +pv -sy) = 0$

(38)  $\alpha (-uq + yt +ys -up) = 0$

(39)  $\alpha (qz- tu +pz -su) = 0$

(40)  $\alpha (vt- zq +vs -zp) = 0$

(41)  $\beta (-pt + vx +ux -qt) = 0$

(42)  $\beta (-xv + sp - xu +sq) = 0$

(43)  $\beta (vs-pz +us -qz) = 0$

(44)  $\beta (-tv + zp +zq -tu) = 0$

(45)  $sq -ux + tq- vx = 0$

(46)  $-up + sy -vp +ty = 0$

(47)  $qu-ys + qv  -yt = 0$

(48)  $-ps + xu-pt  +xv = 0.$

(49)  $ \alpha ut - \alpha zq -px +xq +\alpha pz- \alpha sv = 0$

(50)  $- \alpha vq + \alpha yt  -\beta xt + \beta sx - \alpha sy
+\alpha pu = 0$

(51)  $\beta tp - \beta xv + \alpha yv - \alpha uy - \beta qs+ \beta ux = 0$

(52)  $ - \beta sv + \beta zp+ qy - yp - \beta qz+ \beta ut = 0$

(53)  $pu - ys -\beta tz + \beta zs + ty - vq = 0$

(54)  $-qs + xu + \alpha vz -\alpha zu + tp - vx  = 0$
                   
It is clear that $r$ is the solution of CYBE iff relations (28)-(54) hold.

(I)  By computation, we have that if $r$ is  $\alpha , \beta$-symmetric then
relations (28)-(54) hold and so $r$ is a solution of CYBE.

Conversely, if  $r$ is a solution of CYBE,
then, by computation, we have

(55)  $q^2 =  p^2$   (by  (34)$+$ (32));

(56)  $u^2 =  v^2$  ( by  (33)$+$ (36));

(57)  $t^2 =  s^2$  ( by  (31)$+$ (35) );
                          
Consequently, $p=q, s=t$ and $u=v$ since char  $k =2$.

By relation (31), we have that

$$\alpha y z + \beta xz +xy +   \beta s^2 + \alpha u^2 + p^2 =0.$$
Thus $r$ is  $\alpha ,\beta $-symmetric.

(II)  It is clear that if $r$ is  strongly symmetric then $r$ is $\alpha ,
\beta  $-symmetric. Consequently, $r$ is a solution of CYBE by part (I).
$\Box$

In particular, Proposition \ref {1.4} implies:

\begin {Example} \label {1.5}
 Let $$  sl(2) := \{ x \mid  x  \hbox{ \ is a }  2 \times 2
\hbox { \ matrix
 with trace zero over \  } k   \}$$
 and
 $$e_1 = \left ( \begin {array} {cc}
 0 & 0\\
 1&0
 \end {array}
 \right ),
 e_2 = \left ( \begin {array} {cc}
 0 & 1\\
 0&0
 \end {array}
 \right ),
 e_3 = \left ( \begin {array} {cc}
 1 & 0\\
 0&1
 \end {array}
 \right ).$$
 Thus  $L$ is a Lie algebra (defined by $[x,y] = xy -yx$)
 with a basis $\{ e_1, e_2, e_3 \}.$ It is clear that
 $$[e_1, e_2] = e_3, [e_2, e_3] = 0e_1, [e_3, e_1] =0 e_2.$$
Consequently,
 $r$ is a solution of CYBE iff $r$  is
$0, 0$- symmetric to the basis $\{ e_1, e_2, e_3 \}.$
\end {Example}

\begin {Proposition}
\label {1.6}
$L$ be a Lie algebra with a basis $\{ e_1, e_2, e_3 \}$
such that
$[e_1,e_2]=0,  [e_1,e_3]= e_1+ \beta e_2, [e_2,e_3]=\delta e_2$, where
$\beta, \delta \in k$. Let $ p, q, s, t, u, v, x, y, z \in k$.

(I)  If  $r$ is strongly symmetric,  then
$r$ is a solution of CYBE.

(II) If $\beta =0,  \delta \not=0,$
then $r$   is a solution of CYBE in $L$ iff

  $r=  p(e_1 \otimes e_2)+ q (e_2 \otimes e_1)
+ s (e_1 \otimes e_3)+ s (e_3 \otimes e_1)
+u (e_2 \otimes e_3) +u (e_3 \otimes e_2)
+x (e_1 \otimes e_1) +y (e_2 \otimes e_2) +   z (e_3 \otimes e_3) $

where $ (\delta +1 )zp =(\delta +1 ) qz =(\delta +1) us,    (\delta + 1)u q=
(\delta +1)up$  and $(\delta +1) ps =(\delta +1) qs$

(III) If  $\beta \not= 0 , \delta=1$,
then $r$   is a solution of CYBE in $L$ iff $r$

  $r=  p(e_1 \otimes e_2)+ q (e_2 \otimes e_1)
+ s (e_1 \otimes e_3)+ s (e_3 \otimes e_1)
+u (e_2 \otimes e_3) +u (e_3 \otimes e_2)
+x (e_1 \otimes e_1) +y (e_2 \otimes e_2)
                      +z (e_3 \otimes e_3) $
                      
with $sp =sq, up =qu , zq =zp $ and $s^2=xz.$

(IV)   If  $ \beta = \delta = 0,$ then
$r$ is a  solution  of CYBE  in $L$ iff

$r=  p(e_1 \otimes e_2)+ q (e_2 \otimes e_1)
+ s (e_1 \otimes e_3)+ s (e_3 \otimes e_1)
+u (e_2 \otimes e_3) +v (e_3 \otimes e_2)
+x (e_1 \otimes e_1) +y (e_2 \otimes e_2)
                      +z (e_3 \otimes e_3) $
                      
with $vs =pz, us = qz, qv =pu$ and $(u+v)x = (p+q)s$.

\end {Proposition}

{ \bf Proof.}
Let   $r = \sum_{i,j= 1 }^{3} k_{ij}(e_{i} \otimes e_{j}) \in L \otimes L$
and  $k_{ij} \in k,$ with  $i, j =1, 2, 3.$
 By computation, for all $i, j, n = 1, 2, 3,$   we have that the coefficient of
  $e_j \otimes e_i \otimes e_i$ in
$[r^{12}, r^{13}] $
is zero and
 $e_i \otimes e_i \otimes e_j$
in  $[r^{13}, r^{23}]$ is zero.

We now see the coefficient
 of  $e_i \otimes e_j \otimes e_n$

   in  $[r^{12}, r^{13}] + [r^{12}, r^{23}] +[r^{13}, r^{23}]$.

(1) $e_1 \otimes e_1 \otimes e_1 (
 -k_{13}k_{11} +  k_{11} k_{31})$;
                                  
(2) $e_2 \otimes e_2 \otimes e_2 (
 -\beta k_{23}k_{12} + \beta k_{21} k_{32}
  - \delta k_{23}k_{22} + \delta k_{22} k_{32})$;

(3) $e_3 \otimes e_3 \otimes e_3 ( 0);$

(4) $e_1 \otimes e_2 \otimes e_3 (
 -k_{32}k_{13} +  k_{12} k_{33} -  \beta k_{13}  k_{13}+  \beta k_{11} k_{33}
- \delta k_{13}  k_{23} + \delta k_{12} k_{33} );$

(5) $e_2 \otimes e_3 \otimes e_1 (
 -\beta k_{33}k_{11} + \beta k_{13} k_{31}
 - \delta k_{33}  k_{21}+ \delta k_{23} k_{31}
 -  k_{23}  k_{31} + k_{21} k_{33} );$

(6) $e_3 \otimes e_1 \otimes e_2 (
 -k_{33}k_{12} +  k_{31} k_{32} -  \beta k_{33}  k_{11}+  \beta k_{31} k_{13}
- \delta k_{33}  k_{12} + \delta k_{32} k_{13} );$

(7) $e_1 \otimes e_3 \otimes e_2 (
 -k_{33}k_{12} +  k_{13} k_{32} -  \beta k_{13}  k_{32}+  \beta k_{11} k_{33}
- \delta k_{13}  k_{32} + \delta k_{12} k_{33} );$

(8) $e_3 \otimes e_2 \otimes e_1 (
 -\beta k_{33}k_{11} + \beta k_{31} k_{31} -  \delta k_{33}  k_{21}
 +  \delta k_{32} k_{31}
- k_{33}  k_{21} +  k_{31} k_{23} );$

(9) $e_2 \otimes e_1 \otimes e_3 (
 -\beta k_{31}k_{13} + \beta  k_{11} k_{33}
 - \delta k_{31}  k_{23}+  \delta k_{21} k_{33}
-  k_{23}  k_{13} +  k_{21} k_{33} );$

(10) $e_1 \otimes e_1 \otimes e_2 (
 -k_{31}k_{12} +  k_{11} k_{32} -   k_{13}  k_{12}+   k_{11} k_{32});$

(11) $e_2 \otimes e_1 \otimes e_1 (
 -k_{23}k_{11} +  k_{21} k_{31} -   k_{23}  k_{11}+  k_{21} k_{13});$

(12) $e_1 \otimes e_1 \otimes e_3 (
 -k_{31}k_{13} +  k_{11} k_{33} -   k_{13}  k_{13}+   k_{11} k_{33}
 );$

(13) $e_3 \otimes e_1 \otimes e_1 (
 -k_{33}k_{11} +  k_{31} k_{31} -   k_{33}  k_{11}+  k_{31} k_{13}
);$

(14) $e_2 \otimes e_2 \otimes e_1 (
 -\beta k_{32}k_{11} +  \beta k_{12} k_{31} -  \delta k_{32}  k_{21}
 +  \delta k_{22} k_{31}
- \beta k_{23}  k_{11} + \beta k_{21} k_{31}
- \delta k_{23}  k_{21} + \delta k_{22} k_{31}
 );$

(15) $e_1 \otimes e_2 \otimes e_2 (
 -\beta k_{13}k_{12} + \beta  k_{11} k_{32} -  \delta k_{13}  k_{22}
 +  \delta k_{12} k_{32}
- \beta k_{13}  k_{21} + \beta k_{11} k_{23}
- \delta k_{13}  k_{22} + \delta k_{12} k_{23}
 );$

(16) $e_2 \otimes e_2 \otimes e_3 (
 -\beta k_{32}k_{13} +  \beta k_{12} k_{33} -  \delta k_{32}  k_{23}
 + \delta k_{22} k_{33}
- \beta k_{23}  k_{13} + \beta k_{21} k_{33}
- \delta k_{23}  k_{23} + \delta k_{22} k_{33}
);$

(17) $e_3 \otimes e_2 \otimes e_2 (
 -\beta k_{33}k_{12} + \beta  k_{31} k_{32}- \delta k_{33}k_{22}
 +\delta k_{32}k_{32}
  -  \beta k_{33} k_{21} + \beta k_{31}  k_{23}
  - \delta k_{33} k_{22}
  + \delta k_{32} k_{23}
  );$

(18) $e_3 \otimes e_3 \otimes e_1 (0);$

(19) $e_1 \otimes e_3 \otimes e_3 (0);$

(20) $e_3 \otimes e_3 \otimes e_2 (0);$

(21) $e_2 \otimes e_3 \otimes e_3 (0);$

(22) $e_1 \otimes e_3 \otimes e_1 (
 -k_{33}k_{11} +  k_{13} k_{31} -   k_{13}  k_{31}+   k_{11} k_{33})
 = e_1 \otimes e_3 \otimes e_1 (0)$;

(23) $e_1 \otimes e_2 \otimes e_1 (
 -k_{32}k_{11} +  k_{12} k_{31} -  \beta k_{13}  k_{11}+  \beta k_{11} k_{31}
- \delta k_{13}  k_{21} + \delta k_{12} k_{31}
-  k_{13}  k_{21} +  k_{11} k_{23});$

(24) $e_2 \otimes e_1 \otimes e_2 (
 -\beta k_{31}k_{12} + \beta k_{11} k_{32}
 - \delta k_{31}  k_{22}+  \delta k_{21} k_{32}
- k_{23}  k_{12} +  k_{21} k_{32}
 -\beta k_{23}k_{11} + \beta k_{21} k_{13}
  - \delta k_{23}  k_{12}+  \delta k_{22} k_{13}  );$

(25) $e_2 \otimes e_3 \otimes e_2 (
 -\beta k_{33}k_{12} + \beta k_{13} k_{32}
  -  \delta k_{33}  k_{22}+  \delta k_{23} k_{32}
- \beta k_{23}  k_{31} + \beta k_{21} k_{33}
 + \delta k_{33}  k_{22}-  \delta k_{22} k_{33}
              );$

(26) $e_3 \otimes e_2 \otimes e_3 (
 -\beta k_{33}k_{13} + \beta k_{31} k_{33}
  - \delta k_{33}  k_{23}+  \delta k_{32} k_{33}
);$

(27) $e_3 \otimes e_1 \otimes e_3 (
 - k_{33}k_{13} +   k_{31} k_{33}
 );$

Let $k_{11}=x, k_{22}=y,  k_{33}=z, k_{12}=p,  k_{21} =q, k_{13}=s,
 k_{31} =t, k_{23}=u,  k_{32}=v.$

   It follows from (1)-(27) that

(28) $-sx + xt=0$;

(29) $-\beta up +\beta qv -\delta uy + \delta yv=0$;

(30) $- vs + pz -\beta s^2 + \beta xz -\delta su + \delta zp =0$;

(31) $-\beta xz + \beta st - \delta zq + \delta ut - ut   + qz =0$;

(32) $- zp + tv -\beta zx + \beta ts -\delta zp + \delta vs =0$;

(33) $- zp + sv -\beta st + \beta xz -\delta sv + \delta zp =0$;

(34) $- \beta zx + \beta tt - \delta zq + \delta vt - zq + tu =0$;

(35) $- \beta st + \beta xz  -\delta tu + \delta qz -us + qz=0 $;

(36) $- tp + xv - sp + vx =0$;

(37) $- ux + qt-ux + qs =0$;

(38) $- st + xz - s^2 +  xz =0$;

(39) $-zx +tt -zx +st =0$;

(40) $-\beta vx +\beta pt -\delta vq + \delta yt -\beta ux + \beta qt
-\delta uq + \delta yt =0$;

(41) $-\beta sp + \beta xv  -\delta  sy + \delta pv -\beta sq + \beta xu
- \delta sy + \delta pu =0$;

(42) $-\beta vs + \beta pz -\delta vu + \delta yz -\beta su + \beta zq
-\delta u^2 + \delta yz  =0$;

(43) $-\beta  zp + \beta tv - \delta zy + \delta vv  -\beta zq + \beta tu
-\delta zy +\delta vu =0$;

(44) $- vx + pt -\beta sx + \beta xt -\delta sq + \delta pt
-sq + xu =0$;

(45) $- \beta pt + \beta vx -\delta ty + \delta qv  -up + qv
-\beta ux + \beta qs - \delta up + \delta ys =0$;

(46) $-\beta zp + \beta sv -\beta ut + \beta qz =0$;

(47) $- \beta zs + \beta tz -\delta zu + \delta vz =0$;

(48) $- zs + tz =0$;

It is clear that   $r$   is a solution of CYBE iff
 (28)-(48) hold.

By computation we have that

(49)  $t =s$  (by  (39) + (38) )

 (50)  $\delta u = \delta v$  (by (42) + (43) ).

(I) It is trivial.

(II) Let $\beta =0 $ and $\delta \not=0.$ If $r$ is the case in part (II),
then $r$ is a solution of CYBE by straightforward verification.
Conversely, if $r$ is a solution of CYBE, then $u=v$  by (50) and we have that
\begin {eqnarray*}
(1 +\delta )zp &=& (1+\delta ) us    \hbox { \ \ \ \ by }  (30)   \\
(1 +\delta )zq &=& (1+\delta ) us    \hbox { \ \ \ \ by }  (31)   \\
(1 +\delta )sp &=& (1+\delta ) qs    \hbox { \ \ \ \ by }  (44)   \\
(1 +\delta )qu &=& (1+\delta ) up    \hbox { \ \ \ \ by }  (45)
\end {eqnarray*}
Thus $r$ is the case in part (II).

 (III) Let $\beta \not=0, \delta =1$. If $r$ is the case in part (III)
 then $r$ is a solution of CYBE by straightforward verification. Conversely,
 if $r$ is a solution of CYBE then $u=v$ by (50) and we have that
\begin {eqnarray*}
up &=&qu    \hbox { \ \ \ \ by }  (29)   \\
s^2 &=&xz    \hbox { \ \ \ \ by }  (30)   \\
sp &=&  qs    \hbox { \ \ \ \ by }  (40)   \\
zp &=& zq    \hbox { \ \ \ \ by }  (46)
\end {eqnarray*}
Thus $r$ is the case in part (III).

(IV) Let $\beta =\delta =0.$
It is clear that
the system of equations  (28)-(48) is equivalent to the below
$$  \left    \{  \begin{array} {l}
us =qz   \\
vs =zp  \\
up =qv  \\
(p+q)s = (u+v )x
\end{array} \right.$$
                 Thus we complete the proof.  $\Box$

\begin {Corollary}\label {1.7}
Let $L$ is a Lie algebra with $dim { \ } L \le 3$  and $r \in L \otimes L.$
If $r$ is strongly symmetric, then $r$ is a solution of CYBE in $L$.
                                     
\end{Corollary}

{\bf Proof .}   If $k$ is algebraically closed, then  $r$ is a solution of
 CYBE by Proposition  \ref {1.4}, \ref {1.6} and
 \ref {1.3}.
 If $k$ is not algebraically closed, let
  $P$ be algebraically closure of $k$.
  We can construct a Lie algebra   $ L_P = P \otimes L$ over $P$,
  as in \cite [section 8] {J}.
Set $ \Psi : L \longrightarrow  L_P $  by sending $x $  to
  $1 \otimes x$. It is clear that $L_P$ is a Lie algebra over
  $P$ and  $\Psi$  is  homomorphic with $ker \Psi =0$ over $k.$
  Let  $$\bar r = (\Psi \otimes \Psi )(r).$$
  Obviously, $\bar r$  is  strongly symmetric. Therefore  $\bar r$
 is a solution of CYBE in $L_P$ and so is $r$  in $L$.
$\Box$

\section {Coboundary Lie bialgebras }

In this section, using the general solution, which are obtaied in the section
above,  of CYBE in Lie algebra $L$ with $dim \ \ L \le 3, $ we give the
       the sufficient and  necessary conditions which  $(L, \hbox {[ \ ]},
 \Delta _r, r)$ is a coboundary
 (or triangular ) Lie bialgebra over field $k$ of characteristic 2.
Throughtout this section, the characteristic of field $k$  is 2.

We now observe the  connection between solutions
of CYBE and triangular Lie bialgebra structures.

\begin {Lemma}  \label {2.1.1}
If $r \in Im (1-\tau )$ then $$x\cdot C(r) = (1+\xi +\xi ^2)
 (1 \otimes \Delta) \Delta (x)$$ for any $x \in L$
\end {Lemma}

{\bf Proof:}
Let $r= \sum _i (a_i \otimes b_i - b_i \otimes a_i )$

We see that
\begin {eqnarray*}
C(r) &=& \sum _{i,k}  (
[a_k , a_i] \otimes b_k \otimes b_i
-[a_k , b_i] \otimes b_k \otimes a_i
-[b_k , a_i] \otimes a_k \otimes b_i \\
 &+&  [b_k, b_i] \otimes a_k \otimes a_i
 +  a _k \otimes [b_k, a_i] \otimes b_i
  -a _k \otimes [b_k, b_i] \otimes a_i  \\
  &-&  b _k \otimes [a_k, a_i] \otimes b_i
   -  b _k \otimes [a_k, a_i] \otimes b_i
   +  b _k \otimes [a_k, b_i] \otimes a_i \\
&+&   a _k \otimes a_i \otimes [b_k, b_i]
-  a _k \otimes b_i \otimes [b_k, a_i]
-   b _k \otimes a_i \otimes [a_k, b_i] \\
 &+& b _k \otimes b_i \otimes [a_k, a_i])
 \end {eqnarray*}

and

\begin {eqnarray*}
&{}& ( 1\otimes \Delta )\Delta (x) \\
 &=&
\sum _{i,k} (
[x,a_k] \otimes [b_k,a_i] \otimes b_i +
 [x,a_k] \otimes  a_i \otimes [b_k, b_i]  -
                 [x,a_k] \otimes [b_k,b_i] \otimes a_i \\
&-& [x,a_k] \otimes  b_i \otimes [b_k, a_i]
-             [x,b_k] \otimes [a_k,a_i] \otimes b_i
- [x,b_k] \otimes  a_i \otimes [a_k, b_i] \\
&+&                 [x,b_k] \otimes [a_k,b_i] \otimes a_i
+ [x,b_k] \otimes  b_i \otimes [a_k, a_i]
+a_k \otimes [[x,b_k],a_i] \otimes b_i \\
&+& a_k \otimes  a_i \otimes [[x,b_k], b_i]
 -                    a_k \otimes [[x,b_k],b_i] \otimes a_i
 - a_k \otimes  b_i \otimes [[x,b_k], a_i]   \\
                 &-&              b_k \otimes [[x,a_k],a_i] \otimes b_i
                 -b_k \otimes  a_i \otimes [[x,a_k], b_i]  +
                 b_k \otimes [[x,a_k],b_i] \otimes a_i \\
                 &+&
 b_k \otimes  b_i \otimes [[x,a_k], a_i] )
\end {eqnarray*}

It is easy to check that the sum of the terms whose third factor includes
element $x$ in $(1+\xi +\xi ^2)(1+\Delta )\Delta (x)$ is equal to

\begin {eqnarray*}
 &{}& \sum _{i,k} \{
 [b_k,a_i] \otimes b_i  \otimes [x,a_k] +
 a_i \otimes [b_k, b_i]  \otimes  [x,a_k] -
                [b_k,b_i] \otimes a_i \otimes [x,a_k]  \\
                &-&  b_i \otimes [b_k, a_i]  \otimes [x,a_k]
-             [a_k,a_i] \otimes b_i \otimes  [x,b_k]
-  a_i \otimes [a_k, b_i] \otimes [x,b_k] \\
&+&              [a_k,b_i] \otimes a_i  \otimes  [x,b_k]
+ b_i \otimes [a_k, a_i]  \otimes   [x,b_k]   \\
&+& (b_i \otimes a_k  \otimes [[x,b_k],a_i]
- b_k \otimes  a_i \otimes [[x,a_k], b_i])  \\
&+& ( a_k \otimes  a_i \otimes [[x,b_k], b_i]
 -    a_i \otimes a_k\otimes [[x,b_k],b_i] ) \\
 &+&  ( -  a_k \otimes  b_i \otimes [[x,b_k], a_i]
 +  a_i   \otimes  b_k \otimes  [[x,a_k],b_i] ) \\
&+& (- b_i \otimes  b_k \otimes  [[x,a_k],a_i]
+ b_k \otimes  b_i \otimes [[x,a_k], a_i] ) \}  \\
              &=& (1 \otimes 1 \otimes L_x ) C(r)    \hbox { \ \ \ \
              by Jacobi identity }
\end {eqnarray*}
 where    $L_x$  denotes the adjoint
              action  $L_x(y) = [x,y]$  of $L$ on $L$.

Consequently,   $$(1 +\xi + \xi ^2) (1 \otimes \Delta )\Delta (x)  =
 (L_x \otimes L_x \otimes L_x ) C(r) = x \cdot C(r)$$
for any $x \in L.$  $\Box$

By Lemma \ref {2.1.1} and \cite [ Proposition 2.11 ] {M},
we have:

\begin {Theorem} \label {2.2.1}
$(L, [ \hbox { \ }], \Delta _r ,r)$
is a triangular Lie bialgebra   iff
 $r$ is a solution of CYBE in $L$ and $r \in Im (1-\tau )$
\end {Theorem}

 Consequently, we can easily get a  triangular Lie bialgebra
 structure  by means of a solution of CYBE.

\begin {Theorem}
\label {2.1}
Let $L$ be a  Lie algebra with
       a basis $\{ e_1, e_2, e_3 \}$
such that
$[e_1,e_2]= e_3, [e_2,e_3]=\alpha e_1, [e_3,e_1]=\beta e_2$, where
$\alpha, \beta \in k$.
Then

(I)   $(L, [\hbox { \ }], \Delta _r, r )$ is  a coboundary Lie bialgebra
iff $r \in Im(1 -\tau )$ ;

(II)       $(L, [\hbox { \ }],\Delta _r, r )$ is  a
triangular  Lie bialgebra
iff  $r \in Im (1-\tau )$ and
$r$ is  $\alpha , \beta $- symmetric to the basis of  $\{ e_1, e_2, e_3 \}$
\end {Theorem}

 {\bf Proof.}
(I) Obviously, $r\in Im (1-\tau )$ when
 $(L, [\hbox { \ }], \Delta _r, r )$ is  a coboundary Lie bialgebra. Conversely,
 if $r\in Im (1-\tau )$, let
 $r=  p(e_1 \otimes e_2)- p (e_2 \otimes e_1)
+ s (e_1 \otimes e_3)- s (e_3 \otimes e_1)
+ u (e_2 \otimes e_3)- u (e_3 \otimes e_2)$.
It is sufficient to show that
 $$ (1+ \xi + \xi ^2) (1 \otimes \Delta ) \Delta (e_i) =0$$
for $i=1, 2, 3$  by \cite [Proposition 2.11]{M}.
First, by computation, we have that
\begin {eqnarray*}
&{ \ }& (1 \otimes \Delta ) \Delta (e_1) \\
&=&
\{ ( e_1 \otimes e_2\otimes e_3 ) (\beta ps -  \beta ps)
 +(e_3 \otimes e_1\otimes e_2)(\beta ps)
 + (e_2 \otimes e_3\otimes e_1) (- \beta sp) \} \\
&+&\{ ( e_1 \otimes e_3\otimes e_2 ) (-\beta ps + \beta ps)
 +(e_3 \otimes e_2\otimes e_1)( - \beta ps)
 + (e_2 \otimes e_1\otimes e_3) (\beta sp) \} \\
&+&\{ ( e_1 \otimes e_1\otimes e_2 ) (\alpha \beta su )
 +(e_1 \otimes e_2\otimes e_1)( - \beta \alpha su)
 + (e_2 \otimes e_1\otimes e_1) (0) \} \\
&+& \{ ( e_1 \otimes e_1\otimes e_3 ) ( - \alpha pu)
 +(e_1 \otimes e_3\otimes e_1)(\alpha  pu)
 + (e_3 \otimes e_1\otimes e_1) (0) \} \\
&+&\{   ( e_2 \otimes e_1\otimes e_2 ) (- \beta ^2 s^2)
 +(e_1 \otimes e_2\otimes e_2)(0)
 + (e_2 \otimes e_2\otimes e_1) ( \beta ^2 s^2) \} \\
&+&\{  ( e_3 \otimes e_3\otimes e_1 ) ( p^2 )
 +(e_3 \otimes e_1\otimes e_3)(- p^2)
 + (e_1 \otimes e_3\otimes e_3) (0) \}
\end {eqnarray*}
Thus
 $$ (1+ \xi + \xi ^2) (1 \otimes \Delta ) \Delta (e_1) =0$$
Similarly, we have that
 $$ (1+ \xi + \xi ^2) (1 \otimes \Delta ) \Delta (e_2) =0$$
           $$ (1+ \xi + \xi ^2) (1 \otimes \Delta ) \Delta (e_3) =0$$
Thus  $(L, [\hbox { \ }],\Delta _r,r)$ is a coboundary Lie bialgebra.

(II)
It follows from Theorem \ref {2.2.1} and Proposition \ref {1.4} $\Box$

\begin
{Example} \label {2.2} Under Example \ref {1.5}, we have the
following:
                                                                         
(i)   $(sl(2), [\hbox { \ }],\Delta _r, r )$ is  a coboundary Lie bialgebra
iff $r\in Im (1-\tau )$ ;

(ii)       $(sl(2), [\hbox { \ }],\Delta _r, r )$ is  a
triangular  Lie bialgebra iff
$r$ is  $0 , 0 $- symmetric to the basis of  $\{ e_1, e_2, e_3 \}$
iff
 $r=
 s (e_1 \otimes e_3)+ s (e_3 \otimes e_1)
+ u (e_2 \otimes e_3)+ u (e_3 \otimes e_2)$.
\end {Example}

 \begin {Theorem} \label  {2.3}
Let  $L$ be a Lie algebra
 with a basis   $\{ e_1, e_2, e_3 \}$  such that
$$[e_1,e_2]= 0, [e_1, e_3]= e_1 + \beta e_2, [e_2, e_3]=\delta e_2,$$
where $\delta, \beta \in k$ and $\delta =1$ when $\beta \not=0.$
Let $p, s, u \in k$ and
 $r=  p(e_1 \otimes e_2)- p (e_2 \otimes e_1)
+ s (e_1 \otimes e_3)- s (e_3 \otimes e_1)
+ u (e_2 \otimes e_3)- u (e_3 \otimes e_2)$. Then

(I)  $(L, [\mbox { \ }],\Delta _r, r )$ is a coboundary Lie bialgebra
iff $r\in Im (1-\tau )$ and  $$ (\delta +1) ((\delta +1)u + \beta s)s =0; $$

(II)  $(L, [\mbox { \ }],\Delta _r, r )$ is a triangular Lie bialgebra
 iff $r \in (1-\tau )$ and $$ \beta s +(1+ \delta )  us =0.$$

\end {Theorem}
{ \bf Proof }.
 (I) We get by computation
 \begin {eqnarray*}
&{ \ }& ( 1 \otimes \Delta ) \Delta (e_1) \\
 &=& \{ ( e_1 \otimes e_1\otimes e_2 ) (\beta \delta ss -  \delta us)
 +(e_1 \otimes e_2\otimes e_1)(- \beta \delta ss + \delta us)
 + (e_2 \otimes e_1\otimes e_1) (0) \} \\
&+&  \{ ( e_2 \otimes e_2\otimes e_1 ) (-\beta us + uu+ \beta ^2 s^2
- \beta su)                                       \\
 &+& (e_2 \otimes e_1\otimes e_2)( \beta us -uu + \beta us - \beta ^2 s^2)
  + (e_1 \otimes e_2\otimes e_2) (0) \} \\
&{ \ }& ( 1 \otimes \Delta ) \Delta (e_2) \\
 &=& \{ ( e_1 \otimes e_1\otimes e_2 ) ( \delta ^2 ss)
 +(e_1 \otimes e_2\otimes e_1)(- \delta ^2 s^2 )
 + (e_2 \otimes e_1\otimes e_1) (0) \} \\
 &+& \{ ( e_2 \otimes e_2\otimes e_1 ) ( \delta \beta ss - \delta su)
  +(e_2 \otimes e_1\otimes e_2)( \delta  su - \delta \beta ss )
 + (e_1 \otimes e_2\otimes e_2) (0) \} \\  
&{ \ }& ( 1 \otimes \Delta ) \Delta (e_3)  \\
 &=& \{ ( e_1 \otimes e_2\otimes e_3 ) (\delta su +  \beta ss)
 +(e_2 \otimes e_3\otimes e_1)(-\delta us - \beta ss) \\
 &+& (e_3 \otimes e_1\otimes e_2) ( \delta \delta us + \beta \delta ss
 + \beta ss -su) \} \\
&+&\{ ( e_1 \otimes e_3\otimes e_2 ) (-\delta su-  \beta ss)
 +(e_3 \otimes e_2\otimes e_1)( - \beta \delta  ss - \delta \delta us
 -\beta ss + su ) \\
 &+&  (e_2 \otimes e_1\otimes e_3) (\delta us + \beta ss) \} \\
&+& \{ ( e_1 \otimes e_1\otimes e_2 ) (- \delta  sp - \delta \delta ps
+ sp + \delta sp ) \\
 &+& (e_1 \otimes e_2\otimes e_1)( \delta ps + \delta \delta ps -\delta sp
 -sp) +
 (e_2 \otimes e_1\otimes e_1) (0) \} \\
&+& \{
 ( e_1 \otimes e_1\otimes e_3 ) ( ss)
 +(e_1 \otimes e_3\otimes e_1)(-ss)
 + (e_3 \otimes e_1\otimes e_1) (0) \} \\
&+& \{   ( e_2 \otimes e_2\otimes e_1 ) (- \beta  ps + pu - \delta \beta sp
-\delta \delta up - \delta \beta ps + \delta pu - \beta sp - \delta up)  \\
&+& (e_2 \otimes e_1\otimes e_2)( \beta  ps - pu + \delta \beta sp
+\delta \delta up + \delta \beta ps - \delta pu + \beta sp + \delta up) ) \\
&-& (e_1 \otimes e_2\otimes e_2) (0) \} \\
&+& \{ ( e_2 \otimes e_2\otimes e_3 ) ( \delta \delta uu + \delta \beta us
+ \beta \delta us + \beta \beta ss )   \\
 &+& (e_2 \otimes e_3\otimes e_2)(   -\delta \delta uu - \delta \beta us
- \beta \delta us - \beta \beta ss )
                    +        (e_3 \otimes e_2\otimes e_2) (0) \}
\end {eqnarray*}
Consequently,   $$ (1+ \xi + \xi ^2) (1 \otimes \Delta ) \Delta (e_1) =0$$
 $$ (1+ \xi + \xi ^2) (1 \otimes \Delta ) \Delta (e_2) =0$$
 and
 $$ (1+ \xi + \xi ^2) (1 \otimes \Delta ) \Delta (e_3) =0$$
iff
  $$ \delta ^2 us + \delta \beta s^2 + \beta s^2 - us =0.$$
This implies that
    $(L, [\mbox { \ }],\Delta _r,r)$ is a coboundary Lie bialgebra
    iff    $r\in Im (1-\tau )$ and
 $$ (\delta +1)(( \delta +1) u + \beta s )s =0   { \ \ \ \ }$$

(II) It follows from Theorem \ref {2.2.1}  and Proposition  \ref {1.6}
 $\Box$

\begin {Theorem}  \label   {2.4}
If $L$ is a Lie algebra with $dim L =2$ and $r \in L \otimes L$, then
  $(L, [\mbox { \ }],\Delta _r,r)$ is a triangular  Lie bialgebra   iff
  $(L, [\mbox { \ }],\Delta _r,r)$ is a coboundary   Lie bialgebra
  iff   $r\in Im (1-\tau )$
 \end {Theorem}
 { Proof.}
It is an immediate consequence of the main result of \cite {M}.
  $\Box$

\end {document}